\documentclass[11pt,reqno]{amsart}
\usepackage[T1]{fontenc}
\usepackage{lmodern}
\usepackage{amsmath,amssymb,mathtools}
\usepackage{graphicx}
\graphicspath{{figures/}{output/pdf/figures/}}
\usepackage[margin=1.08in]{geometry}
\usepackage{microtype}
\usepackage{tikz}
\usetikzlibrary{arrows.meta,calc}
\usepackage[colorlinks=true,linkcolor=blue!45!black,citecolor=blue!45!black,urlcolor=blue!45!black]{hyperref}
\hypersetup{pdftitle={The fence conjecture for convex projective polygons},pdfauthor={Erxiao Wang}}
\newtheorem{theorem}{Theorem}[section]
\newtheorem{proposition}[theorem]{Proposition}
\newtheorem{lemma}[theorem]{Lemma}
\theoremstyle{definition}

\theoremstyle{remark}
\newtheorem{remark}[theorem]{Remark}
\numberwithin{equation}{section}
\newcommand{\R}{\mathbb R}
\newcommand{\RP}{\mathbb{RP}}
\newcommand{\dd}{\,\mathrm d}
\newcommand{\Int}{\operatorname{int}}
\newcommand{\Om}{\Omega}
\newcommand{\cl}[1]{\overline{#1}}
\definecolor{fenceblue}{RGB}{40,83,130}
\definecolor{armred}{RGB}{163,56,45}
\definecolor{armgreen}{RGB}{30,109,90}
\tikzset{every picture/.style={line cap=round,line join=round},
  vertex/.style={circle,fill=black,inner sep=1.2pt},
  zero/.style={circle,fill=armred,inner sep=1.7pt}}
\allowdisplaybreaks[1]

\title[The fence conjecture]{The fence conjecture for convex projective polygons}
\author{Erxiao Wang}
\address{School of Mathematical Sciences, Zhejiang Normal University,
Jinhua, Zhejiang Province, China}
\email{wang.eric@zjnu.edu.cn}
\thanks{This work was supported by the National Natural Science Foundation of China (NSFC--RGC grant 12361161603).}
\date{}

\begin{document}
\begin{abstract}
We prove the fence conjecture of Dumas and Wolf: for a convex projective
polygon with at least five vertices, every zero of the Pick differential
lies in the interior of the region obtained by cutting off the triangles
at consecutive triples of vertices. The proof gives a stronger estimate
for arbitrary properly convex domains. In affine-sphere coordinates
normalized at a Pick zero, the domain contains the closed disk of radius
$2^{-1/4}$ and its closure lies in the open disk of radius $2^{1/4}$,
unless it is a rotated square. The proof compares the canonical
Monge--Amp\`ere potential of the domain with that of a normalized square.
Their third-order contact gives alternating local sign sectors. A
planar separation argument and the maximum principle show that the
negative sectors belong to four distinct components, each meeting a
different side of the square. Rotation gives the inner disk bound,
and conormal duality gives the outer bound. Both constants are optimal 
 for the class of properly convex domains normalized at a Pick zero. 
The argument requires neither boundary smoothness nor any restriction
on the multiplicity of the Pick zero.
\end{abstract}
\maketitle

\section{Introduction}

A properly convex domain in $\RP^2$ carries a canonical complete
hyperbolic affine sphere, and hence a Blaschke metric and a holomorphic
cubic differential, the Pick differential. Its zeros are intrinsic
projective invariants of the domain. For polygons, Dumas and Wolf
\cite{DW} identified the corresponding cubic differentials with
polynomial cubic differentials on the complex plane. In particular, an
$n$-gon has $n-3$ Pick zeros counted with multiplicity.

The positions of these zeros are much less explicit than their number.
Dumas and Wolf constructed barriers that keep them away from the
vertices. Their numerical experiments suggested a stronger statement:
the curved barriers could be replaced by the diagonals joining the two
neighbors of each vertex. This is the fence conjecture
\cite[Section 9.3, Conjecture 9.2]{DW}.

Throughout the paper a convex projective $n$-gon means the open interior
$P$ of a properly convex polygon with $n$ distinct extreme vertices
$v_1,\ldots,v_n$, in cyclic order. Indices are taken modulo $n$. Thus no
listed vertex is a redundant point on an edge. In a bounded affine chart
containing $\cl P$, let $H_i$ be the closed half-plane bounded by
$v_{i-1}v_{i+1}$ which does not contain $v_i$. The \emph{fence} is
\begin{equation}\label{eq:fence}
 F(P)=\cl P\cap\bigcap_{i=1}^n H_i .
\end{equation}
Equivalently, one removes the open cap at each consecutive triple of
vertices and retains its opposite side. This construction is independent
of the bounded affine chart, as is membership in the two-dimensional
interior of $F(P)$.

\begin{theorem}[Fence theorem]\label{thm:fence}
Let $P$ be a convex projective polygon with at least five vertices.
Every zero of its Pick differential lies in $\Int F(P)$.
\end{theorem}

The conclusion is strict, including at multiple zeros. As already
observed by Dumas and Wolf \cite[Section 9.3]{DW}, a quadrilateral
has a unique Pick-zero point at its diagonal intersection, and its
fence consists of that point. Triangles have no Pick zeros.

\begin{figure}[tbp]
\centering
\includegraphics[width=.94\textwidth]{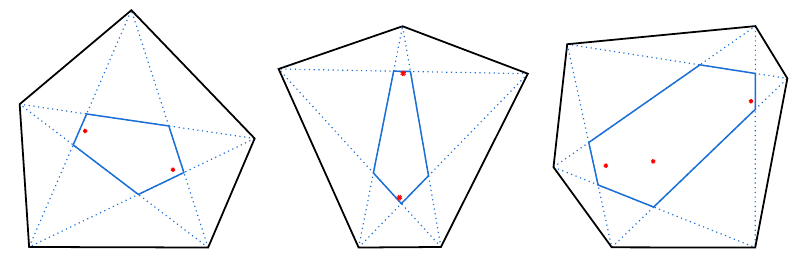}
\caption{Two pentagons and a hexagon from the numerical experiments
of Dumas and Wolf \cite[Figure 4]{DW}. Red points mark the computed
Pick zeros; the inner solid blue polygons are the fences.
The numerical positions are reproduced from their original vector
figure. Theorem~\ref{thm:fence} proves strict containment for every
polygon with at least five vertices; Section~\ref{sec:examples}
discusses these examples and the further families in \cite{DWweb}.}
\label{fig:fences}
\end{figure}

Figure~\ref{fig:fences} reproduces three of the numerical configurations
that motivated the conjecture. Further examples and the role of these
computations are discussed in Section~\ref{sec:examples}.

The estimate underlying Theorem~\ref{thm:fence} applies to domains that
need not be polygonal. We normalize the affine mean curvature to $-1$.
At a point $p$ of the affine sphere, choose a positively oriented
Blaschke-orthonormal tangent frame. The two tangent vectors and the
position vector form a unimodular basis of $\R^3$. Send this basis to
the standard basis and use the affine chart $x_3=1$. We call the
resulting coordinates \emph{coordinates normalized at $p$}. They are
unique up to a rotation. The chart is bounded on the projective domain;
this fact and the corresponding normalization of the potential are
proved in Section~\ref{sec:normalization}. Write $B_\rho$ for the open
Euclidean disk of radius $\rho$ centered at the origin.

\begin{theorem}[Disk bounds at a Pick zero]\label{thm:disks}
Let $\Om\subset\RP^2$ be properly convex, and let $p$ be a zero of the
Pick differential of its complete hyperbolic affine sphere. In
coordinates normalized at $p$, put
\[
 a=2^{-1/4},\qquad R=2^{1/4}.
\]
Exactly one of the following alternatives holds:
\begin{enumerate}
\item $\Om$ is a rotation of the square $(-a,a)^2$;
\item $\cl{B_a}\subset\Om$ and $\cl\Om\subset B_R$.
\end{enumerate}
In the second alternative both inclusions have a positive margin:
there is $\varepsilon>0$, depending on the normalized domain, such that
\[
 \cl{B_{a+\varepsilon}}\subset\Om,\qquad
 \cl\Om\subset B_{R-\varepsilon}.
\]
\end{theorem}

The square attains both radii, so the constants cannot be improved for
the full class of domains in this theorem. We make no assertion of
optimality in a class of polygons with a prescribed number of vertices.

Theorem~\ref{thm:disks} concerns the geometry of an arbitrary properly
convex domain at a point where its cubic form vanishes. Related estimates
compare the Hilbert and Blaschke metrics: Benoist and Hulin
\cite[Proposition 3.4]{BH} established their uniform comparability, and Hildebrand
\cite{Hildebrand} obtained optimal inequalities involving Hilbert
distance and the Blaschke length of projective line segments.
Here the vanishing of the cubic form gives a bound on the entire domain
in the affine chart fixed by the position vector and the Blaschke metric
at that point. The prescribed normalization is essential to the values
of the constants. The resulting inclusions are strict for every domain
other than a square, including domains with nonsmooth boundary, and
their ratio yields the fence theorem for polygons.

Dumas and Wolf's barriers use the homogeneous affine sphere over a
triangular cone \cite[Section 7]{DW}. We instead use a square: its
dihedral symmetry forces the cubic form to vanish at the center,
making it a suitable comparison domain at a Pick zero.
The square model also has a radially symmetric Blaschke metric in a
conformal coordinate in which its Pick differential is $z\,\dd z^3$.
Indeed, rotating $z$ preserves the modulus of this differential, and
uniqueness of the complete, nonpositively curved solution of Wang's
equation \cite[Theorem 5.3]{DW} makes the metric rotation invariant.
This symmetry concerns the conformal coordinate, not the affine
coordinates of the square or its Monge--Amp\`ere potential.
Only the dihedral symmetry and conormal duality are needed below.

The normalized square has the same
third-order potential jet as $\Om$ at $p$. Unless the domains coincide,
the first nonzero term of the difference is a harmonic homogeneous
polynomial of degree at least four. Its alternating signs give at least
four local negative sectors. The essential issue is global: two such
sectors might a priori join elsewhere in the intersection of the
domains. We exclude this by a Jordan-curve argument and the maximum
principle. Each resulting negative component must meet the boundary
of the square at a strictly negative value, while one square side can
serve at most one component. All four sides must therefore meet the
open domain $\Om$. Rotating the square gives the inner disk bound.
Applying the same argument to the conormal dual gives the outer bound.
Their ratio is precisely what is needed to control the diagonals
defining the fence.

\section{Canonical potentials and normalization}\label{sec:normalization}

\subsection{The canonical Monge--Amp\`ere equation}

For a bounded open convex domain $\Om\subset\R^2$, the Cheng--Yau
theorem gives a unique strictly convex function
$u=u_\Om\in C^\infty(\Om)\cap C^0(\cl\Om)$ satisfying
\begin{equation}\label{eq:MA}
 \det D^2u=(-u)^{-4},\qquad u<0\text{ in }\Om,\qquad
 u=0\text{ on }\partial\Om .
\end{equation}
The radial graph
\begin{equation}\label{eq:radial}
 r(x)=-\frac{(x_1,x_2,1)}{u(x)}
\end{equation}
is the complete hyperbolic affine sphere of affine mean curvature $-1$
asymptotic to the cone over $\Om$. We use \eqref{eq:MA} on arbitrary
bounded convex domains, without a boundary smoothness assumption;
see \cite{CY} and the formulation in \cite[Theorem 2.5]{BH}.

Only interior analyticity will be needed. Indeed,
\[
 \log\det D^2u+4\log(-u)=0
\]
is an analytic nonlinear elliptic equation on the set where $u<0$ and
$D^2u$ is positive definite. Its principal linearization is
$(D^2u)^{-1}$, uniformly positive definite on every compact subset of
$\Om$. The analytic interior regularity theorem for nonlinear elliptic
equations therefore shows that $u$ is real analytic in $\Om$; see
\cite[Section 5.8, Theorem 5.8.6]{Morrey}. The application is local at
each interior point, so the degeneracy of \eqref{eq:MA} at
$\partial\Om$ is irrelevant here.

The uniqueness in \eqref{eq:MA} implies invariance under orthogonal
symmetries of $\Om$. Direct substitution also gives the scaling law
\begin{equation}\label{eq:scaling}
 u_{s\Om}(sx)=s^{2/3}u_\Om(x),\qquad s>0.
\end{equation}

\subsection{The jet at a Pick zero}

Let $g$ and $\nabla$ be the Blaschke metric and connection, respectively,
and let $\widehat\nabla$ be the Levi-Civita connection of $g$. We use the
real cubic tensor
\[
 C(X,Y,Z)=g\bigl((\nabla-\widehat\nabla)_X Y,Z\bigr).
\]
Its vanishing is equivalent to the vanishing of the holomorphic Pick
differential. A constant factor in the convention for that differential
does not affect any argument below.

\begin{lemma}\label{lem:jet}
In coordinates normalized at a Pick zero, the projective domain is
bounded, contains $0$, and its canonical potential satisfies
\begin{equation}\label{eq:jet}
 u(0)=-1,\qquad Du(0)=0,\qquad D^2u(0)=I,\qquad D^3u(0)=0.
\end{equation}
\end{lemma}

\begin{proof}
Let $r_1,r_2$ be a positively oriented $g$-orthonormal tangent frame
at the point under consideration. Since the affine normal is $r$,
the Blaschke volume normalization gives
$\det(r_1,r_2,r)=1$. A transformation in $\mathrm{SL}(3,\R)$ sends
these three vectors to $e_1,e_2,e_3$.

The conormal covector $\nu$, characterized by
\[
 \nu(r)=1,\qquad \nu(r_1)=\nu(r_2)=0,
\]
is then $e_3^*$ at the point. Conormal duality identifies the dual
sphere with the complete sphere inside the dual proper cone
\cite[Theorem 3]{Loftin}. Thus $e_3^*$ lies in the interior of that dual
cone: it is strictly positive on every nonzero vector in the closure
of the original cone. The section $x_3=1$ is consequently bounded.
For example, its failure to be bounded would, after division by
Euclidean length and passage to a subsequence, produce a nonzero
closed-cone direction annihilated by $e_3^*$. In this section the
projected point is $0$.

Write $X=(x_1,x_2,1)$ and $r=-X/u$. The position and tangent plane
at $0$ give $u(0)=-1$ and $Du(0)=0$. Differentiation yields
\begin{equation}\label{eq:structure}
 r_{ij}=-\frac{u_i}{u}r_j-\frac{u_j}{u}r_i-\frac{u_{ij}}u r.
\end{equation}
Hence
\begin{equation}\label{eq:gGamma}
 g_{ij}=-\frac{u_{ij}}u,\qquad
 \Gamma^k_{ij}=-\frac{u_i\delta_j^k+u_j\delta_i^k}{u}.
\end{equation}
At $0$ the coordinate tangent vectors are $e_1,e_2$, so that
$g(0)=I$ and $D^2u(0)=I$. Furthermore,
\[
 \Gamma^k_{ij}(0)=0,\qquad
 \partial_k g_{ij}(0)=u_{ijk}(0),\qquad
 \widehat\Gamma^k_{ij}(0)=\tfrac12u_{ijk}(0).
\]
It follows that $C_{ijk}(0)=-\tfrac12u_{ijk}(0)$. The Pick-zero
assumption gives the last equality in \eqref{eq:jet}.
\end{proof}

\subsection{Conormal duality in the normalized chart}

For a bounded convex domain containing $0$, define its open polar by
\begin{equation}\label{eq:polar}
 \Om^\circ=\{y\in\R^2:\ \sup_{x\in\cl\Om}x\cdot y<1\}.
\end{equation}
It is again bounded, convex, and contains $0$. With this convention,
$(\Om^\circ)^\circ=\Om$.

\begin{lemma}\label{lem:dual}
If $u_\Om$ satisfies \eqref{eq:jet}, then the canonical potential of
$\Om^\circ$ satisfies the same normalization at $0$.
\end{lemma}

\begin{proof}
Set $\ell(x)=x\cdot Du(x)-u(x)$. Convexity and $u(0)=-1$
give $\ell(x)\geq1$. Solving $\nu(r)=1$ and $\nu(r_i)=0$
for the radial graph \eqref{eq:radial} gives the conormal covector
\[
 \nu=(-Du,\ell).
\]
Use the unimodular dual basis $(-e_1^*,-e_2^*,e_3^*)$.
A covector with projective coordinates $(y_1,y_2,1)$ then
evaluates on $(x_1,x_2,1)$ as $1-x\cdot y$, so the section of
the dual cone is precisely $\Om^\circ$. The conormal duality
theorem \cite[Theorem 3]{Loftin} identifies $\nu$ with its
canonical affine sphere. Its projective coordinates and
canonical potential are consequently
\begin{equation}\label{eq:dualpotential}
 y=\frac{Du(x)}{\ell(x)},\qquad
 u^\circ(y)=-\frac1{\ell(x)}.
\end{equation}
These formulas also make the normalization explicit. From
\eqref{eq:jet},
\[
 u(x)=-1+\tfrac12|x|^2+O(|x|^4),\qquad
 \ell(x)=1+\tfrac12|x|^2+O(|x|^4),\qquad
 y=x+O(|x|^3).
\]
The last map is locally invertible at $0$. Substituting its
inverse in \eqref{eq:dualpotential} yields
\[
 u^\circ(y)=-1+\tfrac12|y|^2+O(|y|^4),
\]
which is exactly the required third-order jet.
\end{proof}

\begin{remark}[First-order duality]\label{rem:lelieuvre}
The conormal map also has a first-order description. In positively
oriented conformal coordinates $(\sigma,\tau)$ with
$g=e^\psi(\dd\sigma^2+\dd\tau^2)$, it satisfies
\begin{equation}\label{eq:lelieuvre}
 \begin{aligned}
 r_\sigma&=\nu\times\nu_\tau,&
 \nu_\sigma&=r\times r_\tau,\\
 r_\tau&=\nu_\sigma\times\nu,&
 \nu_\tau&=r_\sigma\times r.
 \end{aligned}
\end{equation}
Here the cross products use the ambient volume form and its dual.
These are the Lelieuvre relations of Bobenko--Schief
\cite[Theorem 2.9]{BS}, with signs adjusted to our Gauss equation:
they follow from
$\nu_i(r_j)=-g_{ij}$ and
$\det(r_\sigma,r_\tau,r)=e^\psi$; each right-hand side has
the opposite sign to that in \cite[equation (2--18)]{BS}.
The projective-coordinate formula \eqref{eq:dualpotential} gives
the normalized jet directly. The comparison argument below uses
that formula and the scalar maximum principle.
\end{remark}

\section{The normalized square}\label{sec:square}

Let $S_a=(-a,a)^2$ and let $v_a=u_{S_a}$. The dihedral symmetries
of the square and uniqueness in \eqref{eq:MA} imply
\[
 Dv_a(0)=0,\qquad D^2v_a(0)=\lambda I\quad(\lambda>0),
 \qquad D^3v_a(0)=0.
\]
By \eqref{eq:scaling}, there is exactly one $a>0$ for which
$v_a(0)=-1$. At that value, \eqref{eq:MA} gives $\lambda^2=1$,
and hence $\lambda=1$. This square has the same third-order jet
\eqref{eq:jet} as every domain normalized at a Pick zero.

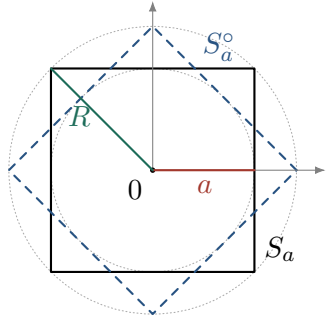
\begin{figure}[tbp]
\centering
\begin{tikzpicture}[scale=1.6]
 \pgfmathsetmacro{\aa}{pow(2,-.25)}
 \pgfmathsetmacro{\rr}{pow(2,.25)}
 \draw[gray!65,densely dotted] (0,0) circle (\aa);
 \draw[gray!65,densely dotted] (0,0) circle (\rr);
 \draw[thick] (-\aa,-\aa) rectangle (\aa,\aa);
 \draw[fenceblue,thick,dashed] (\rr,0)--(0,\rr)--(-\rr,0)--(0,-\rr)--cycle;
 \draw[gray,-{Latex[length=1.4mm]}] (0,0)--(1.43,0);
 \draw[gray,-{Latex[length=1.4mm]}] (0,0)--(0,1.40);
 \fill (0,0) circle (.022);
 \node[below left] at (0,0) {$0$};
 \draw[armred,thick] (0,0)--(\aa,0);
 \node[armred,below] at (.43,0) {$a$};
 \draw[armgreen,thick] (0,0)--(-\aa,\aa);
 \node[armgreen,left] at (-.43,.45) {$R$};
 \node[right] at (\aa,-.65) {$S_a$};
 \node[fenceblue,right] at (.33,1.02) {$S_a^\circ$};
\end{tikzpicture}
\caption{The normalized square and its polar, drawn in the paired
coordinates of Lemma~\ref{lem:dual}. Polarity rotates this square by
$\pi/4$. The dotted circles have radii $a=2^{-1/4}$ and
$R=2^{1/4}$ from \eqref{eq:constants}; this is the equality model
in Theorem~\ref{thm:disks}.}\label{fig:square}
\end{figure}

The size of the square can be determined without an explicit formula
for $v_a$. Its polar is
\[
 S_a^\circ=\{y:\ a(|y_1|+|y_2|)<1\},
\]
a square rotated through $\pi/4$ and having inradius
$a'=1/(\sqrt2\,a)$. By Lemma~\ref{lem:dual}, its canonical potential
also has value $-1$ at the center. The uniqueness of the normalized
size, together with invariance under rotations, forces $a'=a$.
Consequently
\begin{equation}\label{eq:constants}
 a=2^{-1/4},\qquad R=\sqrt2\,a=2^{1/4}=\frac1a .
\end{equation}
Every rotation of $S_a$ therefore has a potential with the jet
\eqref{eq:jet}. Figure~\ref{fig:square} shows the square, its polar,
and the two circles determined by \eqref{eq:constants}.

\section{Square comparison and global separation}\label{sec:comparison}

Fix a domain $\Om$ with potential normalized by \eqref{eq:jet}, and
let $S$ be any rotation of the square in Section~\ref{sec:square}.
Write $u=u_\Om$, $v=u_S$, and
\[
 D=\Om\cap S,\qquad w=u-v.
\]
The common domain $D$ is bounded, open, and convex, and contains $0$.
The function $w$ is analytic in $D$ and continuous on $\cl D$.

\subsection{Contact and the linearized equation}

\begin{lemma}\label{lem:contact}
Either $\Om=S$, or the first nonzero homogeneous Taylor term of $w$
at $0$ is a harmonic polynomial of degree $m\geq4$. In the latter
case, there are $m$ positive and $m$ negative radial segments issuing
from $0$, on which $w$ has alternating strict signs away from $0$.
\end{lemma}

\begin{proof}
If $w$ vanishes in a neighborhood of $0$, analyticity and connectedness
imply $u=v$ on $D$. To see that this forces $\Om=S$, compare the
radial boundary distances of the two domains on any ray from $0$.
If they differ, the first boundary point is interior to the other
domain. Equality of the potentials along the ray would then have
limits $0$ and a strictly negative number, a contradiction.

Suppose $\Om\ne S$. Put $u_t=v+tw$ for $0\leq t\leq1$.
The functions $u_t$ are negative and have positive definite Hessians
on $D$. Subtracting the logarithmic equations for $u$ and $v$ gives
\begin{equation}\label{eq:linear}
 A^{ij}w_{ij}-cw=0,\qquad
 A^{ij}=\int_0^1\bigl(D^2u_t\bigr)^{-1}_{ij}\dd t,\qquad
 c=-4\int_0^1\frac{\dd t}{u_t}>0 .
\end{equation}
Here and below repeated indices are summed from $1$ to $2$.
The matrix $A$ is positive definite and locally uniformly elliptic
in $D$, and $A(0)=I$.

The common jet \eqref{eq:jet} gives $w(x)=O(|x|^4)$.
Since $w$ is analytic and not identically zero, it has a first
nonzero homogeneous term $P_m$, with $m\geq4$. The terms of
degree $m-2$ in \eqref{eq:linear} give $\Delta P_m=0$.
Thus, in polar coordinates,
\begin{equation}\label{eq:harmonic}
 P_m(r,\theta)=b r^m\cos(m\theta-\beta),\qquad b\ne0 .
\end{equation}
Choose one ray strictly inside each of its $2m$ sign sectors.
The expansion $w=P_m+O(|x|^{m+1})$ gives the asserted signs
on sufficiently short segments of these rays.
\end{proof}

We call these segments \emph{local arms}. The number of local arms
does not by itself count the connected components of $\{w>0\}$ or
$\{w<0\}$ on $D$. The following separation argument supplies that
step.

\subsection{Why distinct arms cannot reconnect}

\begin{lemma}[Global separation]\label{lem:separation}
Under the second alternative of Lemma~\ref{lem:contact}, the $m$
positive local arms lie in $m$ distinct components of $\{w>0\}$.
The analogous assertion holds for the negative arms.
\end{lemma}

\begin{figure}[tbp]
\centering
\begin{tikzpicture}[scale=1.25,font=\small]
\begin{scope}[xshift=-2.20cm]
 \fill[armred!8] (0,0)--(1.30,0)
  ..controls(1.55,.85) and (.85,1.55)..(0,1.30)--cycle;
 \draw[gray!50,dashed] (0,0) circle (.49);
 \foreach \t in {135,180,225,270,315}{
  \draw[gray!55,densely dotted] (0,0)--(\t:.68);
 }
 \draw[armred,thick] (0,0)--(1.30,0);
 \draw[armred,thick] (0,0)--(0,1.30);
 \draw[armred,thick] (1.30,0)
   ..controls(1.55,.85) and (.85,1.55)..(0,1.30);
 \draw[armgreen,very thick,-{Latex[length=1.6mm]}] (0,0)--(.57,.57);
 \fill (0,0) circle (.022);
 \node[below left] at (0,0) {$0$};
 \node[armgreen,right] at (.55,.38) {$w<0$};
 \node[armred,right] at (1.30,.82) {$\Gamma$};
 \node[armred,above] at (.45,1.28) {$w>0$};
 \node at (.52,.91) {$G$};
 \node at (.45,-.95) {(a) Reconnection};
\end{scope}
\begin{scope}[xshift=2.15cm]
 \fill[armgreen!9] (-1.05,1.00) rectangle (1.05,1.40);
 \draw[thick] (-1.25,-.50)--(-1.25,1.40)--(1.25,1.40)--(1.25,-.50);
 \draw[gray,dashed] (0,0)--(-1,1.40);
 \draw[gray,dashed] (0,0)--(1,1.40);
 \draw[armgreen,very thick] (-1,1.40)--(1,1.40);
 \draw[armgreen,very thick] (-.75,1.05)--(.75,1.05);
 \node[above] at (0,1.43) {$K\subset E\cap\Om$};
 \node[armgreen,fill=white,inner sep=1pt] at (0,.81) {$(1-\epsilon)K$};
 \node[armgreen] at (0,.40) {$w<0$};
 \fill (0,0) circle (.022);
 \node[below] at (0,0) {$0$};
 \node[right] at (1.26,1.39) {$E$};
 \node at (0,-.95) {(b) One side};
\end{scope}
\end{tikzpicture}
\caption{Two parts of the comparison argument. In (a), the proposed
reconnection of positive arms in Lemma~\ref{lem:separation}
would enclose a negative value with nonnegative boundary
values, contradicting \eqref{eq:linear}. In (b), as in
Lemma~\ref{lem:side}, a compact segment of
$E\cap\Om$ can be moved slightly toward $0$ while retaining $w<0$.
The figures are schematic.}\label{fig:comparison}
\end{figure}
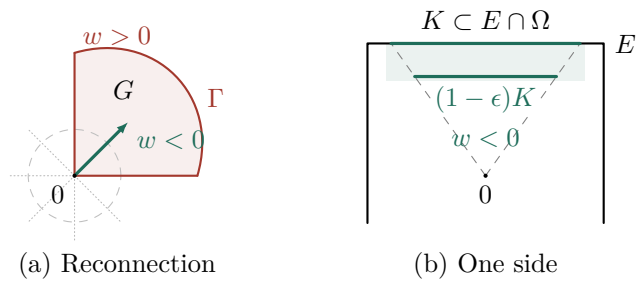

\begin{proof}
Suppose that two distinct positive arms lie in the same positive
component. An open connected subset of the plane is polygonally
path connected, so points on these arms can be joined by a polygonal
path in $\{w>0\}$. Its image is compact and avoids $0$. Removing
loops and making small perturbations within the positive set, we
may arrange that the path is simple and has only finitely many
intersections with the two chosen radial segments.

From this path extract an arc joining the two segments whose interior
meets neither of them: traverse the path to its first encounter with
the second segment, and retain the portion after the last preceding
encounter with the first. The extracted arc, together with the two
radial segments running from its endpoints to $0$, forms a simple
Jordan curve $\Gamma$; see Figure~\ref{fig:comparison}(a).
On $\Gamma\setminus\{0\}$ the function $w$
is positive, while $w(0)=0$.

Let $G$ be the bounded region enclosed by $\Gamma$. Both $\cl G$
and $\Gamma$ lie compactly in $D$. Indeed, the convex hull of the
compact curve is a compact subset of the open convex set $D$,
and contains $\cl G$. Thus the equation \eqref{eq:linear} is valid
throughout a neighborhood of $\cl G$.

Near $0$, the region $G$ occupies one of the two angular regions
between the chosen positive rays. Each of these angular regions
contains a negative local arm: there is a negative ray between any
two successive positive rays in the cyclic list of
Lemma~\ref{lem:contact}. Since the connecting arc stays a positive
distance from $0$, a sufficiently short portion of that negative
arm lies inside $G$. Consequently $w$ is negative somewhere in
$G$ and nonnegative on $\partial G$.

The continuous function $w$ attains a negative minimum at an
interior point of $G$. Its Hessian is positive semidefinite there.
As $A$ is positive definite and $c>0$, this gives
\[
 A^{ij}w_{ij}-cw>0,
\]
contrary to \eqref{eq:linear}. Hence two positive arms cannot
reconnect. Interchanging the signs gives an interior positive
maximum and the opposite strict inequality, excluding the
reconnection of two negative arms as well.
\end{proof}

The proof uses the convexity of the actual intersection domain
$D$, not a larger ambient neighborhood. It also requires no global
description of the nodal set: other zeros and other sign components
may be present. Only the two arms being joined, the connecting
arc, and one intervening arm enter the argument.

\subsection{Strict exits and the capacity of one side}

We next pass from components in the interior to the boundary.
The distinction between a strict sign on the boundary and a
common-boundary zero will be essential.

\begin{lemma}[Strict boundary exits]\label{lem:exits}
Every negative component of $w$ has a point in its closure at which
$w<0$ on $\partial S\cap\Om$. Every positive component has a point
in its closure at which $w>0$ on $\partial\Om\cap S$.
\end{lemma}

\begin{proof}
Let $V$ be a negative component. Its closure in $\cl D$ is compact,
and $w$ has a strictly negative minimum on that closure. If a
minimizing point were in $D$, continuity would give a small ball
about it on which $w<0$. Since this ball meets $V$, it belongs to
the same negative component. The minimum would therefore be an
interior minimum, contradicting \eqref{eq:linear} exactly as in
Lemma~\ref{lem:separation}.

The minimum lies on $\partial D$. The boundary values of the two
potentials give
\begin{equation}\label{eq:boundarysign}
 \begin{aligned}
 w&=-v\geq0 &&\text{on }\partial\Om\cap\cl S,\\
 w&=u\leq0 &&\text{on }\partial S\cap\cl\Om.
 \end{aligned}
\end{equation}
At a point of $\partial\Om\cap\partial S$ both potentials vanish.
A strictly negative boundary minimum must therefore lie on
$\partial S$ and strictly inside $\Om$. The proof for a positive
component uses its positive maximum and the first line of
\eqref{eq:boundarysign}.
\end{proof}

\begin{lemma}[One component per side]\label{lem:side}
Let $E$ be a closed side of $S$. All points of $E\cap\Om$, if any,
are adjacent to the same negative component of $w$ in $D$.
\end{lemma}

\begin{proof}
For $y\in E\cap\Om$ we have $w(y)=u(y)<0$. Continuity gives a
ball about $y$ on whose intersection with $\cl D$ the function
is negative. The intersection with $D$ is nonempty and convex.
Thus $y$ is adjacent to a unique negative component. This
argument includes the case in which $y$ is a vertex of $S$.

Take $y,z\in E\cap\Om$ and let $K=[y,z]$. Convexity shows that
$K$ is a compact subset of $\Om$, and $w=u$ is uniformly
negative on $K$. For all sufficiently small $\epsilon>0$,
\[
 K_\epsilon=(1-\epsilon)K
\]
is a segment in $D$ on which $w<0$. Here moving toward the
interior point $0$ puts $K_\epsilon$ strictly inside $S$;
its containment in $\Om$ and its negative sign follow from
compactness and continuity. Its endpoints lie in the negative
neighborhoods of $y$ and $z$. The segment joins these
neighborhoods within a single negative component.
This inward displacement is illustrated in Figure~\ref{fig:comparison}(b).
\end{proof}

\begin{proposition}[Square comparison]\label{prop:square}
If $\Om\ne S$, every open side of $S$ meets $\Om$, and no vertex
of $S$ lies in $\Om$. Moreover, $w$ has exactly four negative
components and its leading term at $0$ has degree four.
\end{proposition}

\begin{proof}
Let $N_-$ be the number of negative components. Lemmas
\ref{lem:contact} and \ref{lem:separation} give $N_-\geq m\geq4$.
By Lemma~\ref{lem:exits}, every component has a strict exit
through a side of the square, and Lemma~\ref{lem:side} says that
a side can serve at most one component. Thus the entire count is
\[
 4\leq m\leq N_-\leq4.
\]
It follows that $N_-=m=4$ and that all four sides have strict exits,
serving four different components.

If a vertex of $S$ belonged to $\Om$, its negative neighborhood
would identify the components adjacent to the two incident sides.
The four sides could then serve at most three different
components, a contradiction. Thus the exits lie in the relative
interiors of the sides.
\end{proof}

All extrema used above occur either in the interior, where the
operator is elliptic, or on the compact boundary, where only
continuity is used. No boundary Hopf lemma or regularity at a
corner is required.

\section{Disk bounds and the fence}\label{sec:fence}

\subsection{From rotating squares to disks}

\begin{proof}[Proof of Theorem~\ref{thm:disks}]
Suppose first that $\Om$ is not a rotation of $S_a$, with $a$ as
in \eqref{eq:constants}. For $n\in\mathbb S^1$, define the support
function
\[
 h_\Om(n)=\max_{x\in\cl\Om} n\cdot x.
\]
Choose a rotation $S$ of $S_a$ with one side on the line
$n\cdot x=a$. Proposition~\ref{prop:square} gives a point of this
side in the open domain $\Om$. A small displacement in the
direction $n$ stays in $\Om$, so
\begin{equation}\label{eq:support}
 h_\Om(n)>a\qquad(n\in\mathbb S^1).
\end{equation}
The support function is continuous and the circle of directions
is compact. Hence there is $\eta>0$ such that
$h_\Om(n)\geq a+\eta$ for every $n$.

Recall the supporting-half-plane representation
\[
 \cl\Om=\bigcap_{n\in\mathbb S^1}
   \{x:n\cdot x\leq h_\Om(n)\}.
\]
It gives $\cl{B_{a+\eta}}\subset\cl\Om$, and hence
$\cl{B_{a+\eta/2}}\subset\Om$. For the last implication, every
point of $\cl{B_{a+\eta/2}}$ has a ball of radius $\eta/2$
contained in $\cl\Om$, and $\Om=\Int\cl\Om$.

By Lemma~\ref{lem:dual}, $\Om^\circ$ has the same normalization
at $0$. It cannot be a normalized rotated square: the polar of
such a square is another such square, and taking the polar twice
recovers $\Om$. The preceding argument therefore gives
$\cl{B_{a+\delta}}\subset\Om^\circ$ for some $\delta>0$.
For $x\in\cl\Om\setminus\{0\}$, choose
$y=(a+\delta)x/|x|$ in the open polar. Its definition gives
\[
 (a+\delta)|x|=x\cdot y<1.
\]
Thus
\[
 \cl\Om\subset B_{1/(a+\delta)}\subset B_{1/a}=B_R,
\]
with a positive margin relative to $R$. Reducing the two margins
to a common positive $\varepsilon$ proves the second alternative.

If $\Om$ is a normalized rotated square, its side distance from
$0$ is $a$ and its vertices have radius $R$. It satisfies neither
of the strict containments in the second alternative. This also
proves that the alternatives are mutually exclusive.
\end{proof}

\subsection{The angular estimate}

\begin{figure}[tbp]
\centering
\begin{tikzpicture}[scale=1.22,font=\small]
\begin{scope}[xshift=-1.95cm]
 \draw[gray!55,densely dotted] (0,0) circle (1.18921);
 \draw[fenceblue,dashed] (0,0) circle (.84090);
 \draw[thick] (30:1.10)--(90:1.10)--(150:1.10)--
 (210:1.10)--(270:1.10)--(330:1.10)--cycle;
 \draw[armred,very thick] (330:1.10)--(30:1.10);
 \draw[gray,-{Latex[length=1.3mm]}] (0,0)--(1.35,0);
 \node[below] at (1.28,0) {$n_i$};
 \draw (0,0)--(30:1.10);
 \draw (0,0)--(330:1.10);
 \draw[armgreen] (330:.43) arc(-30:30:.43);
 \node[armgreen,fill=white,inner sep=.7pt] at (.43,.02) {$\delta_i$};
 \node[right] at (30:1.10) {$v_{i+1}$};
 \node[right] at (330:1.10) {$v_i$};
 \draw[fenceblue,thin] (0,0)--(135:.8409);
 \node[fenceblue,left] at (-.27,.27) {$a$};
 \node[left] at (155:1.18921) {$R$};
 \node[below left] at (0,0) {$0$};
 \node at (0,-1.50) {(a) $\delta_i<\pi/2$};
\end{scope}
\begin{scope}[xshift=2.0cm]
 \fill[fenceblue!9] (60:1.12)--(120:1.12)--(180:1.12)--
 (240:1.12)--(300:1.12)--cycle;
 \draw[thick] (0:1.12)--(60:1.12)--(120:1.12)--
 (180:1.12)--(240:1.12)--(300:1.12)--cycle;
 \draw[armred,very thick] (300:1.12)--(60:1.12);
 \draw[gray,dashed] (300:1.12)--(0,0)--(60:1.12);
 \draw[gray,dashed] (0,0)--(0:1.12);
 \node[below right] at (300:1.12) {$A$};
 \node[above right] at (60:1.12) {$C$};
 \node[right] at (0:1.12) {$B$};
 \node[zero,label=left:$0$] at (0,0) {};
 \node at (-.34,.42) {$H_i$};
 \node at (0,-1.50) {(b) The retained side of $AC$};
\end{scope}
\end{tikzpicture}
\caption{The passage from disk bounds to the fence. Panel (a)
illustrates the edge estimate of Lemma~\ref{lem:angle} in normalized
coordinates. In (b), as in the proof of Theorem~\ref{thm:fence},
the two consecutive angular gaps from $A$ to $B$ and from $B$ to
$C$ have sum less than $\pi$, placing $0$ opposite $B$ across
the diagonal $AC$.}\label{fig:angles}
\end{figure}
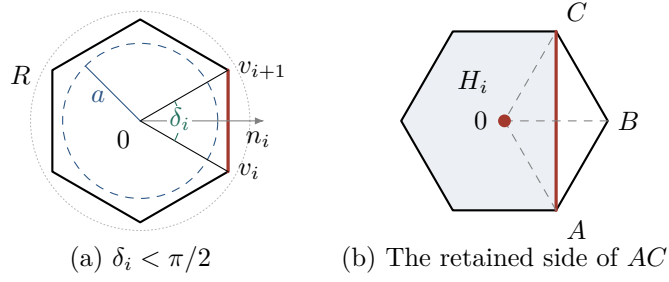

The following elementary consequence explains the particular ratio
of radii in Theorem~\ref{thm:disks}; see Figure~\ref{fig:angles}(a).

\begin{lemma}\label{lem:angle}
Suppose a convex polygon contains $\cl{B_a}$ in its interior and
all its vertices lie in $B_R$, where $a/R=1/\sqrt2$.
Then the positive angular gap between any consecutive vertex
rays from $0$ is strictly less than $\pi/2$.
\end{lemma}

\begin{proof}
Let an edge with endpoints $v_i,v_{i+1}$ have supporting line
$n_i\cdot x=h_i$, where $|n_i|=1$ and $n_i$ points outward.
The inner containment gives $h_i>a$. Both endpoints lie in the
open half-plane $n_i\cdot x>0$, so their signed angles
$\alpha_i,\beta_i$ may be chosen in
$(-\pi/2,\pi/2)$ relative to $n_i$. They satisfy
\begin{equation}\label{eq:angles}
 \cos\alpha_i=\frac{h_i}{|v_i|}>\frac aR=\frac1{\sqrt2},
 \qquad
 \cos\beta_i=\frac{h_i}{|v_{i+1}|}>\frac1{\sqrt2}.
\end{equation}
Thus both angles belong to $(-\pi/4,\pi/4)$. Radial projection
from an interior point sends the convex boundary to the circle
in its cyclic order. In particular, the positive gap along this
edge is the difference of these two angles, not that difference
plus a full turn. It is therefore less than $\pi/2$.
\end{proof}

\begin{proof}[Proof of Theorem~\ref{thm:fence}]
Let $p$ be any Pick zero and normalize at $p$. An invertible
projective transformation preserves the extreme vertices, so a
polygon with at least five vertices cannot become a square.
The strict disk alternative in Theorem~\ref{thm:disks} applies.

Order the vertices counterclockwise, and write $\delta_i$ for the
positive angular gap from $v_i$ to $v_{i+1}$. By
Lemma~\ref{lem:angle},
\[
 0<\delta_i<\frac{\pi}{2}.
\]
Two successive gaps therefore have sum strictly between $0$ and
$\pi$, which gives
\begin{equation}\label{eq:detAC}
 \det(v_{i-1},v_{i+1})>0 .
\end{equation}
Set $A=v_{i-1}$, $B=v_i$, and $C=v_{i+1}$. Relative to the
directed diagonal $AC$, the sign of the origin is
\[
 \det(C-A,-A)=\det(A,C)>0.
\]
Because $A,B,C$ are cyclically ordered extreme vertices of a
convex polygon,
\[
 \det(B-A,C-A)>0,\qquad
 \det(C-A,B-A)<0 .
\]
Thus $0$ lies strictly on the opposite side of $AC$ from $B$,
for every $i$, as illustrated in Figure~\ref{fig:angles}(b).
It is also interior to the polygon. The finitely
many strict inequalities remain valid on a neighborhood of $0$,
so $0\in\Int F(P)$. Returning to the original chart proves the
claim for $p$.\qedhere
\end{proof}

\section{Examples and numerical families}\label{sec:examples}

\subsection{A disk and the equality model}

For the unit disk, the canonical potential is explicit:
\[
 u(x)=-\sqrt{1-|x|^2}.
\]
Its Hessian has tangential eigenvalue $(1-|x|^2)^{-1/2}$ and
radial eigenvalue $(1-|x|^2)^{-3/2}$. Consequently
$\det D^2u=(1-|x|^2)^{-2}=(-u)^{-4}$, and the boundary value
is zero. At the origin its jet is \eqref{eq:jet}; in particular,
\[
 u(x)=-1+\tfrac12|x|^2+\tfrac18|x|^4+O(|x|^6).
\]
The disk lies strictly between the two radii of
Theorem~\ref{thm:disks}, since
$2^{-1/4}<1<2^{1/4}$. Its affine sphere is a quadric, with
identically vanishing cubic form. Thus the disk theorem does
not require the Pick zero to be isolated.

The normalized square realizes the other alternative. Its sides
are tangent to the circle of radius $a$ and its vertices lie on
the circle of radius $R$. Both strict inequalities fail, and the
fence consists of its center. No explicit special-function
formula for its potential is needed to obtain these constants.

\subsection{Regular polygons and numerical examples}

For a regular $n$-gon with vertices $v_j=\rho e_j$, where
$e_j=(\cos(2\pi j/n),\sin(2\pi j/n))$ and $n\geq5$, the fence is
\[
 F(P)=\bigcap_{j=1}^n\{x:x\cdot e_j\leq\rho\cos(2\pi/n)\}.
\]
It is a centered regular polygon of circumradius
$\rho\cos(2\pi/n)/\cos(\pi/n)$. For a pentagon this is
$\rho(3-\sqrt5)/2$. The center is a Pick zero: rotation by
$2\pi/n$ preserves the cubic tensor there and acts on trace-free
cubic tensors by angle $6\pi/n$, with no nonzero fixed vector.

The nonsymmetric examples in Figure~\ref{fig:fences} illustrate
how the Pick zeros sit inside fences of different shapes.
Further numerical families are available in \cite{DWweb}, including
pentagons associated to $(z^2-c)\dd z^3$ and heptagons associated to
$(z^4-c)\dd z^3$. These illustrations are due to Dumas and Wolf;
no new numerical computation is used here. Their website notes
significant numerical error in the older GIF experiments.
Theorem~\ref{thm:fence} includes parameters with repeated roots,
but its strict containment is a statement about each individual
polygon. The domain-dependent margin gives no uniform separation
from the fence along a degenerating family.

\section*{Acknowledgements}

The author thanks J.~F.~Dorfmeister for sustained encouragement in
their work on the Weierstrass representation of affine spheres, even
while the author's research focused on tilings. He is grateful to
the Tsinghua Sanya International
Mathematics Forum for hosting the conference \emph{Surface Theory:
Geometry, Analysis, and Experiments} (25--29 January 2026), which
brought together researchers in Higgs bundles, nonabelian Hodge
theory, differential geometry, and integrable systems.
He thanks H.~Ma, X.~Ma, F.~Pedit, and S.-T.~Yau for the
opportunity to present ongoing work on mixed-type affine spheres. He 
thanks A. Bobenko, S. Heller and L. Heller for helpful discussions on affine spheres. 
He also thanks Q.~Li for introducing him to Dumas and Wolf's
correspondence between convex polygonal cones and polynomial cubic
differentials on $\mathbb C$, and to their fence conjecture, which
inspired this work.

DeepSeek, Doubao, and ChatGPT assisted with proof exploration,
mathematical checks, literature searches, and exposition.
ChatGPT/Codex also assisted with exploratory numerical tests;
these tests are not used as evidence for the results proved here.
The reproduced numerical illustrations are due to Dumas and Wolf.
The author takes full responsibility for the mathematical content.

\enlargethispage{\baselineskip}

\end{document}